\documentclass{amsart}
\usepackage{mathrsfs,pstricks}
\usepackage[colorlinks,
           linkcolor=red,
           anchorcolor=blue,
           citecolor=green]{hyperref}
\usepackage{amsfonts}
\usepackage{enumerate}
\topskip  \usepackage{latexsym}
\usepackage{amsmath,amssymb}
\usepackage{amsthm}
\newtheorem{theorem}{Theorem}[section]

\newtheorem{corollary}[theorem]{Corollary}
\newtheorem{definition}[theorem]{Definition}

\newtheorem{lemma}[theorem]{Lemma}

\newtheorem{property}[theorem]{Property}

\begin{document}
\title{Total embedding distributions of Ringel ladders}
\maketitle
\begin{center}
Yichao Chen, Lu Ou, Qian Zou\\[6pt]
College of mathematics and econometrics, Hunan University, 410082 Changsha, China\\[6pt]

{\it ycchen@hnu.edu.cn, 50371081@qq.com,
Joe\_king520@qq.com}
\end{center}
\footnote{The work was partially supported by NNSFC under
Grant No. 10901048}
\section*{abstract}
The total embedding distributions of a graph is consisted of the orientable embeddings
and non-orientable embeddings and have been know for few classes of graphs.
The genus distribution of Ringel ladders is determined in
[Discrete Mathematics 216 (2000) 235-252] by E.H. Tesar.
In this paper, the explicit formula for  non-orientable embeddings of
 Ringel ladders is obtained.

\medskip

\noindent {\bf Key words}: Graph embedding; Ringel ladders; Overlap matrix; Chebyshev polynomials;

\noindent {\bf 2000 Mathematics Subject Classification}:  05C10, 30B70, 42C05

\section{Background}\label{background}

One enumerative
aspect of topological graph theory is to count genus distributions of a graph.
The history of genus distribution began with J. Gross in 1980s. Since then, it has been attracted a lot of attentions,
 for the details, we may refer to
 \cite{Arc,Ed60,GF87,GKP10,GRT89,Jac87,FGS89,KPG10,KL93,McG87,PKG10,Rie,Sta90,Sta91a,Sta91b,Sta96c,Tes00,VW07, WL06,WL08} etc (We only list a few).
However, for the total embedding distributions, only few
classes are known. For example, Chen, Gross and Rieper \cite{CGR94} computed the
total embedding distribution for necklaces of type $(r,0)$, close-end ladders and cobblestone
paths, Kwak  and Shim \cite{KS02} computed for bouquets of circles and
dipoles. In \cite{CLW06}, Chen, Liu and Wang calculated the total embedding
distributions of all graphs with maximum genus 1.
Furthermore, in  \cite{CMZ09}, Chen, Mansour and Zou obtained
explicit formula for total embedding distributions for
the necklaces of type $(r,s),$ closed-end ladders and cobblestone path.

It is assumed that the reader is somewhat familiar with the basics of
topological graph theory as found in Gross and Tucker \cite{GT87}.
A graph $G=(V(G),E(G))$ is permitted to have both
loops and multiple edges. A \textit{surface} is a compact closed
2-dimensional manifold without boundary. In topology, surfaces are
classified into $O_m$, the \textit{orientable surface} with $m
(m\geq 0)$ handles and $N_n$, the \textit{nonorientable surface}
with $n(n>0)$ crosscaps. A graph embedding into a surface means a
\textit{cellular embedding}.

A \textit{spanning tree} of a graph $G$ is a tree on its edges has the same order as
$G$. The number co-tree edges of a spanning tree of $G$ is called the \textit{Betti number}, $\beta(G)$, of $G$. A
\textit{rotation} at a vertex $v$ of a graph $G$ is a cyclic order
of all edges incident with $v$. A \textit{pure rotation system
}$P$ of a graph $G$ is the collection of rotations at all
vertices of $G$. A \textit{general rotation system }is a pair
$(P,\lambda)$, where $P$ is a pure rotation system and $\lambda$ is
a mapping $E(G)\rightarrow\{0,1\}$. The edge $e$ is said to be {\it twisted}
(respectively, {\em untwisted}) if $\lambda(e)=1$ (respectively, $\lambda(e)=0$). It is well known
that every orientable embedding of a graph $G$ can be described by a
general rotation system $(P,\lambda)$ with $\lambda(e)=0$ for all $e\in E(G)$. By allowing $\lambda$ to take the non-zero
value, we can describe nonorientable embeddings of $G$, see \cite{CGR94,Sta78} for more details.
A \textit{$T$-rotation system }$(P,\lambda)$
of $G$ is a general rotation system $(P,\lambda)$ such that
$\lambda(e)=0$, for all $e \in E(T)$.

\begin{theorem}\label{tha} (see \cite{CGR94,Sta78}) Let $T$ be a
spanning tree of $G$ and $(P,\lambda)$ a general rotation system.
Then there exists a general rotation system $(P^{'},\lambda^{'})$
such that
\begin{itemize}
\item[(1)] $(P^{'},\lambda^{'})$ yields the same embedding of $G$
as $(P,\lambda)$, and
\item[(2)] $\lambda^{'}(e)=0$, for all  $e\in E(T)$.
\end{itemize}
\end{theorem}

Two embeddings are considered to be the
\textit{same} if their $T$-rotation systems are combinatorially equivalent.
Fix a spanning tree $T$ of a graph $G$. Let $\Phi_{G}^{T}$
be the set of all $T$-rotation systems of $G$. It is known that
\[|\Phi_{G}^{T}|=2^{\beta(G)}\prod_{v\in V(G)}(d_{v}-1)!.\]
Suppose that in these $|\Phi_{G}^{T}|$ embeddings of $G$, there are
$a_i$, $i=0,1,\ldots$, embeddings into orientable surface $O_i$
and $b_j$, $j=1,2,\ldots$, embeddings into nonorientable surface
$N_j$. We call the polynomial
\[I_{G}^{T}(x,y)=\sum_{i= 0}^{\infty} a_i x^i +\sum_{j=1}^{\infty} b_j
y^j\] the \textit{$T$-distribution polynomial} of $G$.
By the \textit{total genus polynomial }of $G$, we shall mean the
polynomial
\[I_G(x,y)=\sum_{i= 0}^{\infty} g_i x^i +\sum_{i=1}^{\infty} f_i y^i,\]
where $g_i$ is the number of embeddings (up to equivalence) of $G$
into the orientable surface $O_i$ and $f_i$ is the number of
embeddings (up to equivalence) of $G$ into the nonorientable surface
$N_i$. We call the first (respectively, second) part of $I_G(x,y)$ the \textit{genus polynomial}
(respectively, \textit{crosscap number polynomial}) of $G$ and denoted by
$g_G(x)=\sum_{i= 0}^{\infty} g_i x^i$ (respectively, $f_G(y)=\sum_{i=1}^{\infty} f_i y^i$).
Clearly, $I_G(x,y)=g_G(x)+f_G(y)$. This means the  number of orientable embeddings of
$G$  is $\prod_{v\in G}(d_v-1)!,$ while the number of non-orientable embeddings of $G$ is
$(2^{\beta(G)}-1)\prod_{v\in G}(d_v-1)!.$

Let $T$ be a spanning tree of $G$ and $(P^{'},\lambda^{'})$ be a
$T$-rotation system. Let $e_1,e_2,\ldots,e_{\beta(G)}$ be the cotree
edges of $T$. The \textit{overlap matrix} of $(P^{'},\lambda^{'})$
is the $\beta\times \beta$ matrix $M=[m_{ij}]$ over $GF(2)$ such
that $m_{ij}=1$ if and only if either $i\neq j$ and the restriction
of the underlying pure rotation system to $T+e_i+e_j$ is nonplanar,
or $i=j$ and $e_i$ is twisted. The following theorem due to Mohar.

\begin{theorem}\label{Moh89} (see \cite{Moh89}) Let $(P,\lambda)$ be a general rotation system for a graph, and let $M$
be the overlap matrix. Then the rank of $M$ equals twice the genus,
if the corresponding embedding surface is orientable, and it equals
the crosscap number otherwise. It is independent of the choice of a
spanning tree.
\end{theorem}

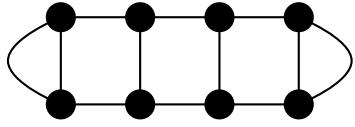
\begin{figure}[h]
\begin{pspicture}(-1.5,-1.6)(4.5,1.)
\psset{xunit=20pt,yunit=22pt}
\psline(0,0)(4.5,0)\pscircle*(0,-1.5){0.2}\pscircle*(1.5,-1.5){0.2}\pscircle*(3,-1.5){0.2}\pscircle*(4.5,-1.5){0.2}
\psline(0,-1.5)(4.5,-1.5)\pscircle*(0,0){0.2}\pscircle*(1.50,0){0.2}\pscircle*(3,0){0.2}\pscircle*(4.5,0){0.2}
\psline(0,-1.5)(0,0) \psline(1.5,-1.5)(1.5,0)\psline(3,-1.5)(3,0)\psline(4.5,0)(4.5,-1.5)
\pscurve(0,0)(-1,-0.75)(0,-1.5)
\pscurve(4.5,0)(5.5,-0.75)(4.5,-1.5)
\end{pspicture}
\caption{The $4$-rung closed-end ladder $L_4$}\label{fig1}
\end{figure}

An {\it $n$-rung closed-end ladder $L_n$} can be obtained by taking
the graphical cartesian product of an {\it $n$-vertex} path with the
complete graph $K_2$, and then doubling both its end edges. Figure
\ref{fig1} presents a $4$-rung closed-end ladder.

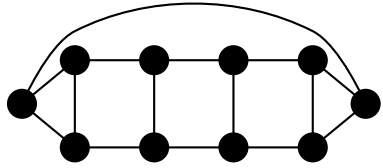
\begin{figure}[h]
\begin{pspicture}(-1.5,-1.6)(4.5,1.)
\psset{xunit=20pt,yunit=22pt}
\psline(0,0)(4.5,0)\pscircle*(0,-1.5){0.2}\pscircle*(1.5,-1.5){0.2}\pscircle*(3,-1.5){0.2}\pscircle*(4.5,-1.5){0.2}
\psline(0,-1.5)(4.5,-1.5)\pscircle*(0,0){0.2}\pscircle*(1.50,0){0.2}\pscircle*(3,0){0.2}\pscircle*(4.5,0){0.2}
\psline(0,-1.5)(0,0) \psline(1.5,-1.5)(1.5,0)\psline(3,-1.5)(3,0)\psline(4.5,0)(4.5,-1.5)
\pscircle*(5.5,-0.75){0.2} \pscircle*(-1,-0.75){0.2}

\psline(0,0)(-1,-.75)\psline(0,-1.50)(-1,-.75) \psline(5.5,-.750)(4.5,0)\psline(5.5,-.750)(4.5,-1.50)

\pscurve(-1,-.75)(0,0.5)(4.5,0.5)(5.5,-.75)
\end{pspicture}
\caption{The Ringel ladder $R_4$}\label{fig2}
\end{figure}

Ringel ladders, $R_n$, are the graphs used by Ringel and Youngs in their proof of the Heawood Map
Coloring Theorem. In fact, A\textit{ Ringel ladder}, $R_n$, can be formed by subdividing the end-rungs of
the closed-end ladder, $L_n$, and adding an edge between these two new vertices. Figure \ref{fig2} shows
the Ringel ladder $R_4$.

\section{Homogeneous recurrence relation and Chebyshev polynomials}\label{Homo}

To begin with the discussion, we give some concepts of
the $n$-th Chebyshev polynomials of the second kind which is related
to the solution of the recurrence relation.
Let the recurrence  function $U_n(x)$ be
$$U_n(x)=2xU_{n-1}(x)-U_{n-2}(x)$$
with the initial conditions $U_0(x)=1$, $U_1(x)=2x$, then we derived
the {\it $n$-th Chebyshev polynomials with the second kind}
$U_n(x)$ (see \cite{Ri}). For instance, $U_2(x)=4x^2-1$, $U_3(x)=8x^3-4x$,
$U_4(x)=16x^4-12x^2+1$. Moreover, we have the identity that
\begin{align}
\label{gs0}
 U_n(x)=\sum\limits_{k=0}^{\lfloor n/2\rfloor}
\binom{n-k}{k}(-1)^{k}(2x)^{n-2k}.
\end{align}
Now, we will build the relation between the recurrence relation
 and the Chebyshev polynomials with the second kind. Let
$P_n(z)=\sum\limits_{m=0}^{n}C_n(m)z^{m},$  satisfy the following
$$P_n(z)=a_1(z)P_{n-1}(z)+a_2(z)P_{n-2}(z),$$
where $a_i(z)=\sum\limits_{k=0}^{q}a_{i,\ k}z^{k}$ for $i=1,2.$ and the initial conditions $P_0(z)=c_{0},$ and $ P_1(z),\ P_2(z)$ can be derived
 by the initial values of $C_{n}(m).$\\
 Let $Q_n(z)=\frac{P_n(z)}{(\sqrt{a_2(z)}i)^{n}},$ then it is easy
 to verify that
 $$Q_n(z)=\frac{a_1(z)}{\sqrt{a_2(z)}i}Q_{n-1}(z)-Q_{n-2}(z)$$
 with the initial conditions
 $Q_0(z)=P_{0}(z)=c_{0},\ Q_1(z)=\frac{P_1(z)}{\sqrt{a_2(z)}i}$ and
 $Q_2(z)=\frac{P_2(z)}{-a_2(z)}.$ Using the fact that $U_0(x)=1,\ U_1(x)=2x,\
 U_2(x)=4x^2-1,$ by induction on $n=0,1,2,$ we obtain that
\begin{align}\label{gs1}
Q_n(z)=AU_n(\frac{a_1(z)}{2\sqrt{a_2(z)}i})+BU_{n-1}(\frac{a_1(z)}{2\sqrt{a_2(z)}i})+CU_{n-2}(\frac{a_1(z)}{2\sqrt{a_2(z)}i}),
\end{align}
where
$A,B,$ and $C$ are determined by the initial conditions.

Thus we have
\begin{align}\label{gs2}
P_n(z)=(\sqrt{a_2(z)}i)^{n}AU_n(\frac{a_1(z)}{2\sqrt{a_2(z)}i})+BU_{n-1}(\frac{a_1(z)}{2\sqrt{a_2(z)}i})+CU_{n-2}(\frac{a_1(z)}{2\sqrt{a_2(z)}i}).
\end{align}
Using the fact that $$U_n(x)=\sum\limits_{k=0}^{\lfloor n/2\rfloor}
\binom{n-k}{k}(-1)^{k}(2x)^{n-2k}.$$ We can derive that
\begin{align}\
(i\sqrt{a_2(z)})^{n}U_n(\frac{a_1(z)}{2\sqrt{a_2(z)}i})=\sum\limits_{j\geq
0}\binom{n-j}{j}(a_1(z))^{n-2j}.
\end{align}
Since $a_1(z)$ is a polynomial of degrees less than $q,$ then
$(a_1(z))^{n-2j}$ can be expressed as the type of power series. Plug
the above formula into (\ref{gs2}) and comparing the coefficient
$z^m$ in both sides and we can obtain the explicit formulae $C_n(m)$
for $0\leq m\leq n$.

\section{Total embedding distributions of Ringel ladders}

\subsection{The rank-distribution polynomial of Closed-end ladders}

we adopt the notations of \cite{CMZ09}, the overlap of matrix of Closed-end ladders $L_{n-1}$ has
 the following form $M_n^{X,Y}$( see \cite{CMZ09} for more details).

Let $X=(x_1,x_2,\ldots,x_n)\in(GF(2))^{n}$ and $Y=(y_1,y_2,\ldots,y_{n-1})\in(GF(2))^{n-1}$. We define the tridiagonal matrix $M_n^{X,Y}$ as
$$M_n^{X,Y}=\left( \begin{array}{lllllll} x_1&y_1& & & & & \\
y_1& x_2&y_2& & &\bf{0}& \\
 &y_2&x_3&y_3& & & \\
 & & & & & & \\
 &\bf{0}& & &y_{n-2}&x_{n-1} &y_{n-1} \\
 & & & & &y_{n-1}&x_n
  \end{array}
  \right).$$
Furthermore, we define
$\mathscr{L}_n=$$\{M_n^{X,Y}\mid X\in (GF(2))^{n}\mbox{ and }Y\in
(GF(2))^{n-1}\}$, which is the set of all matrices over $GF(2)$
that are of the type $M_n^{X,Y}$. We define the {\it rank-distribution polynomial} to be
the polynomial $\mathscr{L}_n(z)=\sum_{j=0}^{n}D_n(j)z^{j}$, where $D_n(j)$, $j=0,1,\ldots,n$, is the number of different assignment
of the variables $x_j,y_k$, where $j=1,2,\cdots,n$ and $k=1,2,\cdots,n-1$, for which
the matrix $M_n^{X,Y}$ in $\mathscr{L}_n$ has rank $j$.
Similarly, Let
$\mathscr{O}_n=\{M_n^{0,Y}\mid Y\in
(GF(2))^{n-1}\}$,  and $\mathscr{O}_n(z)=\sum_{j=0}^{n}O_n(j)z^{j}$ be the {\it rank-distribution polynomial} of $\mathscr{O}_n$
,where $O_n(j)$, $j=0,1,\ldots,n$, is the number of different assignment
of the variables $y_k$, where $k\in \{1,2,\ldots,n-1\}$, for which
the matrix $M_n^{Y}$ in $\mathscr{A}_n$ has rank $j$.

\begin{lemma}\label{added}(see \cite{CMZ09})
The polynomial $\mathscr{O}_n(z)$ satisfies the recurrence relation
$$\mathscr{O}_{n}(z)=\mathscr{O}_{n-1}(z)+2z^{2}\mathscr{O}_{n-2}(z)$$
with the initial conditions $\mathscr{O}_{1}(z)=1$ and $\mathscr{O}_{2}(z)=z^2+1$.
\end{lemma}

\begin{theorem}\label{coa1}(see \cite{CMZ09})
For all $n\geq1$,
\begin{align*}
\mathscr{O}_n(z)&=\sum\limits_{j\geq 0}\binom{n-j}{j}2^j\ z^{2j}-\sum\limits_{j\geq
0}\binom{n-2-j}{j}2^j\ z^{2j+2}.
\end{align*}
\end{theorem}

\begin{corollary}
\label{coro33}For all $1\geq m\leq [\frac{n}{2}].$
\begin{align*}
 &\mathscr{O}_n(2m+1)=0,\\
 &\mathscr{O}_n(2m)=\binom{n-m}{m}\cdot2^m-\binom{n-m-1}{m-1}\cdot 2^{m-1}.
 \end{align*}

 \end{corollary}

\begin{lemma}\label{ttha}(see \cite{CMZ09})
The polynomial $\mathscr{L}_n(z)$ satisfies the recurrence relation
$$\mathscr{L}_{n}(z)=(1+2z)\mathscr{L}_{n-1}(z)+4z^{2}\mathscr{L}_{n-2}(z)$$
with the initial conditions $\mathscr{L}_{1}(z)=1+z$ and $\mathscr{L}_{2}(z)=4z^2+3z+1$.
\end{lemma}

\begin{theorem}\label{coa1}(see \cite{CMZ09})
For all $n\geq1$,
$$\mathscr{L}_n(z)=(2iz)^{n}\left[U_{n}\left(\frac{1+2z}{4iz}\right)+\frac{i}{2}U_{n-1}\left(\frac{1+2z}{4iz}\right)-\frac{1}{2}U_{n-2}\left(\frac{1+2z}{4iz}\right)\right],$$
where $U_s(t)$ is the $s$-th Chebyshev poynomial of the second kind
and $i^2=-1$.
\end{theorem}

\begin{corollary}\label{coa2} (see \cite{CMZ09})
For all $n\geq1$ and $0\leq m\leq n$,
\begin{align*}
D_n(m)&=2^m\sum_{j=0}^{[m/2]}\binom{n-j}{j}\binom{n-2j}{n-m}-2^{m-1}\sum_{j=0}^{[(m-1)/2]}\binom{n-1-j}{j}\binom{n-1-2j}{n-m}\\
&+2^{m-1}\sum_{j=0}^{[(m-2)/2]}\binom{n-2-j}{j}\binom{n-2-2j}{n-m}.
\end{align*}
\end{corollary}

\subsection{The overlap matrix of Ringel ladders}

We adopt the same notation used by Ringel [27, p.17]. A cubic graph at each vertex has two
cyclic orderings of its neighbors. One of these two cyclic orderings is denoted as  clockwise and the
other \textit{counterclockwise}. We color the vertex\textit{ black}, if that vertex has the \textit{clockwise} ordering of its neighbors,
otherwise, we will color the counterclockwise vertices \textit{white}. This will bring
convenient to embed a cubic graph into surfaces, as we can draw an imbedding
on the plane and only need to color the vertices black and white.
\begin{definition}
An edge is called \textit{matched }if it has the same color at both ends, otherwise
it is called \textit{unmatched}.
\end{definition}
We fix a spanning tree $T$ of $R_{n-1}$ shown as the
thicker lines in Figure \ref{fig3}, that is to say, the cotree edges are $e,a_1,a_2,\cdots,a_n$.
\begin{property}\label{p1}
Two cotree edges $e$ and $a_i$, for $i=1,2,\cdots,n$, overlap if and only if the edge $c_i$ is unmatched.
\end{property}

\begin{property}\label{p2}
Two cotree edges $a_i$ and $a_{i+1}$,for $i=1,2,\cdots,n-1$, overlap if and only if the edge $b_i$ is unmatched.
\end{property}

\begin{figure}[h]
\begin{pspicture}(-1.5,-2.2)(8,1.5)

\psline[linewidth=.1cm](0,0)(6,0)\pscircle*(0,-1.5){0.2}\pscircle*(1.5,-1.5){0.2}\pscircle*(3,-1.5){0.2}\pscircle*(4.5,-1.5){0.2}
\psline(0,-1.5)(6,-1.5)\pscircle*(0,0){0.2}\pscircle*(1.50,0){0.2}\pscircle*(3,0){0.2}\pscircle*(4.5,0){0.2}
\psline[linewidth=.1cm](0,-1.5)(0,0) \psline[linewidth=.1cm](1.5,-1.5)(1.5,0)\psline[linewidth=.1cm](3,-1.5)(3,0)\psline[linewidth=.1cm](4.5,0)(4.5,-1.5)
\pscircle*(7.5,-0.75){0.2} \pscircle*(6,0){0.2}\pscircle*(6,-1.5){0.2}
\psline[linewidth=.1cm](0,0)(-1.5,-.75)\psline(0,-1.50)(-1.5,-.75) \psline[linewidth=.1cm](7.5,-.750)(6,0)\psline(7.5,-.750)(6,-1.50)
\psline[linewidth=.1cm](6,0)(6,-1.5)
\pscurve(-1.5,-.75)(-.5,0.5)(6.5,0.5)(7.5,-.75)
\put(3.5,-.75){{$\cdots$}}
\put(-.95,-1.5){$a_1$} \put(0.5,-1.8){$a_2$} \put(2,-1.8){$a_3$} \put(5,-1.8){$a_{n-1}$} \put(6.5,-1.5){$a_{n}$}
\put(-.5,-.75){$b_1$} \put(1.,-.75){$b_2$} \put(2.5,-.75){$b_3$} \put(5.2,-.75){$b_{n-1}$}
\put(-.95,-.15){$c_1$} \put(0.5,.25){$c_2$} \put(2,0.25){$c_3$} \put(5,0.25){$c_{n-1}$} \put(6.5,-.15){$c_{n}$}
\put(3,1.){$e$}
\pscircle*(-1.5,-0.75){0.2}

\put(-.15,.45){$v_2$} \put(1.35,.45){$v_3$} \put(2.85,.45){$v_4$}  \put(4.25,.45){$v_{n-1}$}  \put(5.85,.45){$v_{n}$}
\put(7.5,-.45){$v_{n+1}$} \put(-1.7,-.45){$v_{1}$}

\put(-.15,-2){$u_2$} \put(1.35,-2){$u_3$} \put(2.85,-2){$u_4$}  \put(4.25,-2){$u_{n-1}$}  \put(5.85,-2){$u_{n}$}

\end{pspicture}
\caption{}\label{fig3}
\end{figure}

It is easy to see that the overlap matrix  of  $R_{n-1}$  has the following form.

$$M_{n+1}^{X,Y,Z}=\left( \begin{array}{cccccccc}
x_0&z_1&z_2 & z_3&\ldots &z_{n-1} & z_n& \\
z_1&x_1&y_1& & & & &  \\
z_2&y_1& x_2&y_2& &\bf{0} && \\
z_3& &y_2&x_3& \ddots& & & \\
\vdots & & & \ddots&\ddots & y_{n-2}& & \\
z_{n-1}& &\bf{0}&  &y_{n-2} &x_{n-1} &y_{n-1} \\
z_n& & & & &y_{n-1}&x_n
\end{array}\right),$$
where $X=(x_0,x_1,\ldots,x_n)\in (GF(2))^{n},$
$Y=(y_1,y_2,\ldots,y_{n-1})\in (GF(2))^{n-1}$ and $Y=(z_1,z_2,\ldots,z_{n})\in (GF(2))^{n-1}$.
Note that $x_0=1$ if and only if the edge $e$ is twisted, $x_i=1$ if and only if the edge $a_i$ is twisted, for all
$i=1,2,\ldots,n$, $y_j=1$ if and only if $b_j$ is unmatched.
for all $j=1,2,\ldots,n-1$, and $z_k=1$ if and only if $c_k$ is unmatched,
for all $k=1,2,\ldots,n$.
\begin{property}\label{p3}
For a fixed matrix of the form $M_{n+1}^{X,Y,Z}$, there are exactly
$2$ different $T$-rotation systems
corresponding to that matrix.
\end{property}
\begin{proof}Given a matrix $M_{n+1}^{X,Y,Z}$, the values of $z_1,z_2,\cdots,z_n$ and $y_1,y_2,\cdots,y_{n-1}$ are determined.

\begin{itemize}
  \item $z_1=0.$ If we color the vertex $v_1$ black, by Property \ref{p1}, the color of  $v_2$  is black. Since the values of
  $z_2,\cdots,z_n$ and $y_1,y_2,\cdots,y_{n-1}$ are given, by  Property \ref{p1} and Property \ref{p2}, all the colors of
  $v_2,u_2,\cdots,v_n,u_n,v_{n+1}$ are determined. That is to say, all the rotations of vertices of $R_n$ is determined.
   Otherwise the vertex $v_1$ is colored white, by Property \ref{p1}, the color of
   $v_2$  is also white, by the values of  $z_2,\cdots,z_n$ and $y_1,y_2,\cdots,y_{n-1}$ and by Property \ref{p1} and Property \ref{p2},
   the color all vertices of $R_n$ is determined.
  \item $z_1=1$, Similar discuss like the case $z_1=0$, the details are omitted.

\end{itemize}
\end{proof}

Now, we denote $\mathscr{R}_{n+1}$ be the set of all matrices over $GF(2)$ that are of the form  $M_{n+1}^{X,Y,Z}$.
The we calculate the rank distribution of the set $\mathscr{R}_{n+1}$.

Let $\mathscr{R}_{n+1}(z)=\sum_{j=0}^{n+1}C_{n+1}(j)z^{j}$ be the {\it rank-distribution polynomial} of the set  $\mathscr{R}_{n+1}$.
In other words, for $j=0,1,\ldots,n+1$, $C_{n+1}(j)$ is the number of different assignment
of the variables $x_i$,  $i=0,1,\cdots,n$, $y_k,$ $k=1,2,\cdots,n-1$, and $z_l,$ $l=1,2,\cdots,n$ for
which the matrix $M_{n+1}^{X,Y,Z}$ in  $\mathscr{R}_{n+1}$ has rank $j$.

Similarly, Let $\mathscr{P}_{n+1}$ be the set of all matrices over $GF(2)$ that are of the form  $M_{n+1}^{O,Y,Z}$.
The we calculate the rank distribution of the set $\mathscr{P}_{n+1}$. Let $\mathscr{P}_{n+1}(z)=\sum_{j=0}^{n+1}D_{n+1}(j)z^{j}$ be the {\it rank-distribution polynomial} of the set  $\mathscr{O}_{n+1}$.
In other words, for $j=0,1,\ldots,n+1$, $D_{n+1}(j)$ is the number of different assignment
of the variables  $y_k,$ $k=1,2,\cdots,n-1$, and $z_l,$ $l=1,2,\cdots,n$ for
which the matrix $M_{n+1}^{O,Y,Z}$ in  $\mathscr{P}_{n+1}$ has rank $j$.

\begin{lemma}\label{lemm}
The polynomial $\mathscr{P}_n(z)$ $(n\geq3)$ satisfies the recurrence relation
\begin{align}
\mathscr{P}_{n+1}(z)=\mathscr{P}_{n}(z)+8z^2\mathscr{P}_{n-1}(z)+2^{n-1}z^2\mathscr{O}_{n-1}(z).
\end{align}
with the initial condition  $\mathscr{P}_2(z)=z^2+1,\
\mathscr{P}_{3}(z)=7z^2+1$ and $\mathscr{P}_{4}(z)=12z^4+19z^2+1$ where $\mathscr{O}_{n-1}(z)$ is rank-distribution polynomial of closed-end ladders $L_{n-2}.$
\end{lemma}

\begin{proof}
To obtain the relation between $\mathscr{P}_{n+1}(z)$ and $\mathscr{P}_{n}(z)$, we consider the four different ways to assign the variables
$y_{n-1}$ and $z_n$ in the matrix $M_{n+1}^{Y,Z}$.

\textbf{Case 1:} $y_{n-1}=0$.
\begin{itemize}
  \item

\textbf{Subcase 1:} $z_n=0.$ Then the rank of $M_{n+1}^{Y,Z}$ is the same as the upper left $n\times n$ submatrix, which
is a matrix of the form  $M_n^{Y,Z}.$ We conclude that this case contributes to the polynomial $\mathscr{P}_{n+1}(z)$ by a term $\mathscr{P}_{n}(z).$

\item \textbf{Subcase 2:} $z_n=1.$
 It is easy to sea that, no matter what assignments of the variables $z_{1},z_2,\cdots,$ $z_{n-1},$
we can transform $M_{n+1}^{Y,Z}$ to the following form.
\end{itemize}
$$M_1=\left( \begin{array}{cccccccc}
0&0&0 & 0&\ldots &0 & 1& \\
0&0&y_1& & & & &  \\
0&y_1& 0&y_2& & && \\
0& &y_2&0& \ddots& & & \\
\vdots & & & \ddots&\ddots & y_{n-2}& & \\
0& &&  &y_{n-2} &0 &0 \\
1& & & & &0&0
\end{array}\right),$$

We firstly delete the first column and the last column then delete the first row and the last row of $M_1$, then
we obtain a matrix which is a overlap matrix of closed ladders $L_{n-2}$.
Since there are $2^{n-1}$ different assignments of the variables $z_{1},z_2,\cdots,$ $z_{n-1},$
it contributes to the polynomial $\mathscr{P}_{n+1}(z)$ by
 a term $2^{n-1}z^2\mathscr{O}_{n-1}(z).$

\textbf{Case 2:}$y_{n-1}=1$. If $z_n=1$, we first add the last row to the first low, then add the last column to the fist column.   A similar discussion for $y_{n-2}$ and $z_{n-1}$,  we transform $M_{n+1}^{Y,Z}$ to the
following form.
$$M=\left( \begin{array}{cccccccc}
0&z_1&z_2 & \ldots& z_{n-2} &0 & 0& \\
z_1&0&y_1& & & & &  \\
z_2&y_1& 0&\ddots& & && \\
\vdots & &\ddots&\ddots& y_{n-3}& & & \\
z_{n-2}& & & y_{n-3}&0&0& & \\
0& &&  &0 &0 &1 \\
0& & & & &1&0
\end{array}\right),$$

 Note that the upper left $(n-1)\times (n-1)$ submatrix of $M_2$, which
is a matrix of the form  $M_{n-1}^{Y,Z}.$ There are $2^3$ different assignments of the variables
 $y_{n-2}$, $z_{n-1}$ and $z_n$ in the matrix $M_n^{Y,Z}$. In this case,
it contributes to the polynomial $\mathscr{P}_{n+1}(z)$ by
 a term $8z^2\mathscr{P}_{n}(z).$
\end{proof}

\begin{lemma}\label{lemm}
The polynomial $\mathscr{R}_n(z)$ $(n\geq3)$ satisfies the recurrence relation
\begin{align}
\mathscr{R}_{n+1}(z)=(4z+1)\mathscr{R}_{n}(z)+16z^2\mathscr{R}_{n-1}(z)+2^{n}z^2\mathscr{L}_{n-1}(z).
\end{align}
with the initial condition  $\mathscr{R}_2(z)=4z^2+3z+1,\
\mathscr{R}_{3}(z)=28z^3+28z^2+7z+1,$ where $\mathscr{L}_{n-1}(z)$ is rank-distribution polynomial of closed-end ladders $L_{n-2}.$
\end{lemma}

\begin{proof}
To obtain the relation between $\mathscr{R}_{n+1}(z)$ and $\mathscr{R}_{n}(z)$, we consider the eight different ways to assign the variables
$x_n,$ $y_{n-1}$ and $z_n$ in the matrix $M_{n+1}^{X,Y,Z}$.

\textbf{Case 1:} $x_n=0$.
\begin{itemize}
  \item

\textbf{Subcase 1:} $y_{n-1}=z_n=0.$ Then the rank of $M_{n+1}^{X,Y,Z}$ is the same as the upper left $n\times n$ submatrix, which
is a matrix of the form  $M_n^{X,Y,Z}.$ We conclude that this case contributes to the polynomial $\mathscr{R}_{n+1}(z)$ by a term $\mathscr{R}_{n}(z).$

\item \textbf{Subcase 2:} $y_{n-1}=z_n=1.$
We first add the last row to the first low, then add the last column  to the fist column. If
$x_{n-1}=1,$ we add the last column to the $n$-th column. A similar discussion for $y_{n-2}$ and $z_{n-1}$,  we transform $M_{n+1}^{X,Y,Z}$ to the
following form.
$$M=\left( \begin{array}{cccccccc}
x_0&z_1&z_2 & \ldots& z_{n-2} &0 & 0& \\
z_1&x_1&y_1& & & & &  \\
z_2&y_1& x_2&\ddots& & && \\
\vdots & &\ddots&\ddots& y_{n-3}& & & \\
z_{n-2}& & & y_{n-3}&x_{n-2} &0& & \\
0& &&  &0 &0 &1 \\
0& & & & &1&0
\end{array}\right),$$

Note that the upper left $(n-1)\times (n-1)$ submatrix, which
is a matrix of the form  $M_{n-1}^{X,Y,Z}.$ There are $2^3$ different assignments of the variables
 $x_{n-1},$ $y_{n-2}$ and $z_{n-1}$,  in these case it contributes to the polynomial $\mathscr{R}_{n+1}(z)$ by
 a term $8z^2\mathscr{R}_{n-1}(z).$
 \item
\textbf{Subcase 3:} $y_{n-1}=1, z_n=0.$ Similarly discuss like subcase 2, it  contributes to the polynomial $\mathscr{R}_{n+1}(z)$ by
 a term $8z^2\mathscr{R}_{n-1}(z).$
\item
\textbf{Subcase 4:} $y_{n-1}=0, z_n=1.$ It is easy to sea that, no matter what assignments of the variables $x_{0},$ $z_{1},z_2,\cdots,$ $z_{n-1},$
we can transform $M_{n+1}^{X,Y,Z}$ to the following form.
\end{itemize}
$$M_1=\left( \begin{array}{cccccccc}
0&0&0 & 0&\ldots &0 & 1& \\
0&x_1&y_1& & & & &  \\
0&y_1& x_2&y_2& & && \\
0& &y_2&x_3& \ddots& & & \\
\vdots & & & \ddots&\ddots & y_{n-2}& & \\
0& &&  &y_{n-2} &x_{n-1} &0 \\
1& & & & &0&0
\end{array}\right),$$

We firstly delete the first column and the last column then delete the first row and the last row of $M_1$, then
we obtain a matrix which is a overlap matrix of closed ladders $L_{n-2}$.
Since there are $2^n$ different assignments of the variables $x_{0},$ $z_{1},z_2,\cdots,$ $z_{n-1},$
it contributes to the polynomial $\mathscr{R}_{n+1}(z)$ by
 a term $2^nz^2\mathscr{L}_{n-1}(z).$

\textbf{Case 2:} $x_n=1$. If $z_n=1,$ we first add the last column to the first column then add the last row to the first row.
 Similarly, if  $y_{n-1}=1$, we add the last column to the $n$-th column and add the last row to the $n$-th row. As last we can
 transfer the matrix $M_{n+1}^{X,Y,Z}$ to the matrix $M_2$ of following form.

$$M_2=\left( \begin{array}{cccccccc}
x_0&z_1&z_2 & z_3&\ldots &z_{n-1} & 0& \\
z_1&x_1&y_1& & & & &  \\
z_2&y_1& x_2&y_2& & && \\
z_3& &y_2&x_3& \ddots& & & \\
\vdots & & & \ddots&\ddots & y_{n-2}&0 & \\
z_{n-1}& & &  &y_{n-2} &x_{n-1} &0 \\
0& & & & &0&1
\end{array}\right),$$

 Note that the upper left $n\times n$ submatrix of $M_2$, which
is a matrix of the form  $M_{n}^{X,Y,Z}.$ There are $2^2$ different assignments of the variables
 $y_{n-1}$ and $z_n$ in the matrix $M_n^{X,Y,Z}$. In this case,
it contributes to the polynomial $\mathscr{R}_{n+1}(z)$ by
 a term $4z\mathscr{R}_{n}(z).$
\end{proof}

\begin{theorem}\label{th}
For all $n\geq2$,
$$\aligned
\mathscr{P}_n(z)=&(2\sqrt{2}iz)^{n}\left\{U_n(\frac{1}{4\sqrt{2}iz})-\frac{z^2+1}{4\sqrt{2}iz}
U_{n-1}(\frac{1}{4\sqrt{2}iz})+\frac{17z^2-1}{16z^2}U_{n-2}(\frac{1}{4\sqrt{2}iz})\right\}\\&+2^{n-1}z^2\mathscr{O}_{n-1}(z).
\endaligned
$$
where $U_s(t)$ is the $s$-th Chebyshev poynomial of the second kind, $i^2=-1$ and  $\mathscr{O}_{n-1}(z)$ is rank-distribution polynomial of closed-end ladders $L_{n-2}.$
\end{theorem}
\begin{proof}
Note that
\begin{align}
\label{lem1}
\mathscr{P}_{n+1}(z)=\mathscr{P}_{n}(z)+8z^2\mathscr{P}_{n-1}(z)+2^{n-1}z^2\mathscr{O}_{n-1}(z).
\end{align}

We first consider the homogeneous recurrence relation part of (\ref{lem1}).
\begin{align}
\label{lemt}
\mathscr{P}_{n+1}(z)=\mathscr{P}_{n}(z)+8z^2\mathscr{P}_{n-1}(z).
\end{align}
By the method of subsection \ref{Homo}, we have a solution of (\ref{lemt}).
 \begin{align}
\label{gs6}
 \mathscr{P}_n(z)=(\sqrt{a_2(z)}i)^{n}\left\{AU_n(\frac{a_1(z)}{2\sqrt{a_2(z)}i})+BU_{n-1}(\frac{a_1(z)}{2\sqrt{a_2(z)}i})+CU_{n-2}(\frac{a_1(z)}{2\sqrt{a_2(z)}i})\right\}
\end{align}

Now, let $Y_n(z)=2^{n}f(z)\mathscr{O}_n(z)$ be one special solution of
$\mathscr{P}_n(z),$ plug it into (\ref{lem1}), using the relation $$\mathscr{O}_{n}(z)=(1+2z)\mathscr{O}_{n-1}(z)+4z^{2}\mathscr{O}_{n-2}(z),$$
it leads to
$$Y_n(z)=2^{n-1}z^2\mathscr{O}_{n-1}(z)=\sum\limits_{m\geq 0}2^{n-1}O_{n-1}(m)z^{m+2}.$$

Thus,
\begin{align}
\label{gs8}
\mathscr{P}_n(z)=(2\sqrt{2}zi)^{n}\left\{U_n(\frac{1}{2\sqrt{2}iz})+BU_{n-1}(\frac{1}{2\sqrt{2}iz})+CU_{n-2}(\frac{1}{2\sqrt{2}iz})\right\}+2^{n-1}z^2\mathscr{O}_{n-1}(z).
\end{align}
Plug the initial values $\mathscr{P}_2(z), \mathscr{P}_3(z)$ into (\ref{gs8}), it
follows that
\[ \left\{
\begin{array}{l}
-8z^2\left\{(2(\frac{1}{2\sqrt{2}iz})+B)\ \frac{1}{2\sqrt{2}iz}+(C-1)\right\}+2z^2=z^2+1\\
                                             \\
-16\sqrt{2}iz^3\left\{(\frac{1}{2\sqrt{2}iz}+B)\
(-\frac{1}{8z^2}-1)+\frac{1}{2\sqrt{2}iz})(C-1)\right\}+4z^2(z^2+1)=7z^2+1.
\end{array}
\right.
\]
By simple computation, we immediately obtain
\begin{align*}
&B=\frac{-z^2-1}{4\sqrt{2}iz},\ \
C=\frac{17z^2-1}{16z^2}.
\end{align*}

\end{proof}

Then according to the identity (\ref{gs0}), the formula (\ref{gs8})
is as follows
\begin{align}
\mathscr{P}_n(z)=&\sum\limits_{j\geq
0}\binom{n-j}{j}(8z^2)^{j}-\frac{z^2+1}{2}\times\left\{\sum\limits_{j\geq
0}\binom{n-1-j}{j}(8z^2)^{j}\right\}\notag\\
&-\frac{17z^2-1}{2}\left\{\sum\limits_{j\geq\notag
0}\binom{n-2-j}{j}(8z^2)^{j}\right\}+2^{n-1}z^2\mathscr{O}_{n-1}(z)\notag.
\end{align}
Comparing the coefficient of $z^m$ in both sides,
thus for all $n\geq 2$ and $0\leq m \leq n,$ we have the following result.

\begin{theorem}\label{th}
For all $n\geq2$,
$$\aligned
\mathscr{R}_n(z)=&(4zi)^{n}\left\{U_n(\frac{1+4z}{8iz})-\frac{2z^2+7z+1}{8iz}
U_{n-1}(\frac{1+4z}{8iz})+\frac{34z^2-z-1}{32z^2}U_{n-2}(\frac{1+4z}{8iz})\right\}\\&+2^nz^2\mathscr{L}_{n-1}(z).
\endaligned
$$
where $U_s(t)$ is the $s$-th Chebyshev poynomial of the second kind, $i^2=-1$ and  $\mathscr{L}_{n-1}(z)$ is rank-distribution polynomial of closed-end ladders $L_{n-2}.$
\end{theorem}
\begin{proof}
Note that
\begin{align}
\label{lem1}
\mathscr{R}_n(z)=(4z+1)\mathscr{R}_{n-1}(z)+16z^2\mathscr{R}_{n-2}(z)+2^{n-1}z^2\mathscr{L}_{n-2}(z).
\end{align}

We first consider the homogeneous recurrence relation part of (\ref{lem1}).
\begin{align}
\label{lemt}
\mathscr{R}_n(z)=(4z+1)\mathscr{R}_{n-1}(z)+16z^2\mathscr{R}_{n-2}(z).
\end{align}
By the method of subsection \ref{Homo}, we have a solution of (\ref{lemt}).
 \begin{align}
\label{gs6}
 \mathscr{R}_n(z)=(\sqrt{a_2(z)}i)^{n}\left\{AU_n(\frac{a_1(z)}{2\sqrt{a_2(z)}i})+BU_{n-1}(\frac{a_1(z)}{2\sqrt{a_2(z)}i})+CU_{n-2}(\frac{a_1(z)}{2\sqrt{a_2(z)}i})\right\}
\end{align}

Now, let $Y_n(z)=2^{n}f(z)\mathscr{L}_n(z)$ be one special solution of
$\mathscr{R}_n(z),$ plug it into (\ref{lem1}), using the relation $$\mathscr{L}_{n}(z)=(1+2z)\mathscr{L}_{n-1}(z)+4z^{2}\mathscr{L}_{n-2}(z),$$
it leads to $f(z)=\frac{z^2\mathscr{L}_{n-1}(z)}{\mathscr{L}_n(z)}.$

Thus we obtain a
special solution of non-homogeneous recurrence (\ref{lem1})
$$Y_n(z)=2^nz^2\mathscr{L}_{n-1}(z)=\sum\limits_{m\geq 0}2^{n}C_{n-1}(m)z^{m+2}.$$

Thus,
\begin{align}
\label{gs8}
\mathscr{R}_n(z)=(4zi)^{n}\left\{U_n(\frac{1+4z}{8iz})+BU_{n-1}(\frac{1+4z}{8iz})+CU_{n-2}(\frac{1+4z}{8iz})\right\}+2^nz^2\mathscr{L}_{n-1}(z).
\end{align}
Plug the initial values $\mathscr{R}_2(z), \mathscr{R}_3(z)$ into (\ref{gs8}), it
follows that
\[ \left\{
\begin{array}{l}
-16z^2\left\{(2(\frac{1+4z}{8iz})+B)\ U_1(\frac{1+4z}{8iz})+C-1\right\}+4z^2(1+z)=4z^2+3z+1\\
                                             \\
-64iz^3\left\{(2(\frac{1+4z}{8iz})+B)\
U_2(\frac{1+4z}{8iz})+U_1(\frac{1+4z}{8iz})(C-1)\right\}+8z^2(4z^2+3z+1)\\=28z^3+28z^2+7z+1.
\end{array}
\right.
\]
By simple computation, we immediately obtain
\begin{align*}
&B=-\frac{-2z^2-7z-1}{8iz},\ \
C=\frac{34z^2-z-1}{32z^2}.
\end{align*}

\end{proof}

Then according to the identity (\ref{gs0}), the formula (\ref{gs8})
is as follows
\begin{align}
\label{gs9} \mathscr{R}_n(z)=&\sum\limits_{j\geq
0}\binom{n-j}{j}(1+4z)^{n-2j}(4z)^{2j}-\frac{2z^2+7z+1}{2}\times\\&\left\{\sum\limits_{j\geq
0}\binom{n-1-j}{j}(1+4z)^{n-1-2j}(4z)^{2j}\right\}\notag\\
&+\frac{34z^2-z-1}{2}\left\{\sum\limits_{j\geq\notag
0}\binom{n-2-j}{j}(1+4z)^{n-2-2j}(4z)^{2j}\right\}\\&+2^{n}z^2\mathscr{L}_{n-1}(z)\notag.
\end{align}
Comparing the coefficient of $z^m$ in both sides of (\ref{gs9}),
thus for all $n\geq 2$ and $0\leq m \leq n,$ we have the following result.
\begin{corollary}\label{co}
For all $n\geq2$ and $0\leq m\leq n$,
\begin{align*}
C_n(m)&=\sum\limits_{j=0}^{\lfloor \frac
{m}{2}\rfloor}\binom{n-j}{j}\binom{n-2j}{n-m}4^m
-\sum\limits_{j=0}^{\lfloor \frac
{m-2}{2}\rfloor}\binom{n-j-1}{j}\binom{n-1-2j}{n-m+1}4^{m-2}\\&-\frac{7}{2}\sum\limits_{j=0}^{\lfloor
\frac {m-1}{2}\rfloor}\binom{n-j-1}{j}\binom{n-1-2j}{n-m}4^{m-1}-\frac{1}{2}\sum\limits_{j=0}^{\lfloor \frac
{m}{2}\rfloor}\binom{n-j-1}{j}\binom{n-1-2j}{n-m-1}4^{m}\\
&-17\sum\limits_{j=0}^{\lfloor
\frac {m-2}{2}\rfloor}\binom{n-j-2}{j}\binom{n-2-2j}{n-m}4^{m-2}
+\frac{1}{2}\sum\limits_{j=0}^{\lfloor \frac
{m-1}{2}\rfloor}\binom{n-j-2}{j}\binom{n-2-2j}{n-m-1}4^{m-1}\\&+\frac{1}{2}\sum\limits_{j=0}^{\lfloor
\frac {m}{2}\rfloor}\binom{n-j-2}{j}\binom{n-2-2j}{n-m-2}4^{m}+2^{n}D_{n-1}(m-2).
\end{align*}
where
\begin{align*}
D_n(m)&=2^m\sum_{j=0}^{[m/2]}\binom{n-j}{j}\binom{n-2j}{n-m}-2^{m-1}\sum_{j=0}^{[(m-1)/2]}\binom{n-1-j}{j}\binom{n-1-2j}{n-m}\notag\\
&+2^{m-1}\sum_{j=0}^{[(m-2)/2]}\binom{n-2-j}{j}\binom{n-2-2j}{n-m}.
\end{align*}
\end{corollary}

\begin{theorem} The total genus polynomial of
Ringel ladders $R_{n-1}$ is as follows:
\begin{align*}
\mathbb{I}_{R_{n-1}}(x,y)=2\sum_{j=0}^{n+1}C_{n+1}(j)y^{j} -\mathbb{I}_{0}(R_{n-1},\ y^{2})+\mathbb{I}_{0}(R_{n-1},\ x)
\end{align*}
where $\mathbb{I}_{0}(R_{n-1},x)$ is the
genus polynomial of  Ringel ladder $R_{n-1}$,
which has been derived by E.H. Tesar ~\cite{Tes00}.
\end{theorem}
\begin{proof}
By Property \ref{p3}, the theorem follows.
\end{proof}
For instance, the above theorem gives

\begin{align*}I_{R_1}(x,y)=&2+14x+14y+42y^2+56y^3,\\
I_{R_2}(x,y)=&2+38x+24x^2+22y+122y^2+424y^3+392y^4,\\
I_{R_3}(x,y)=&2+70x+184x^2+30y+242y^2+1448y^3+3272y^4+2944y^5,\\
I_{R_4}(x,y)=&2+118x+648x^2+256x^3\\&+38y+410y^2+3496y^3+12952y^4+26880y^5+20736y^6,\\
I_{R_5}(x,y)=&2+198x+1656x^2+2240x^3\\&+46y+642y^2+7240y^3+36808y^4+120832y^5+207168y^6+147456y^7.
\end{align*}


\end{document}